\documentclass{amsart}
\title{Measurable cardinals and the cardinality of Lindel\"of spaces}
\author{Marion Scheepers}
\date{Version of January 29, 2010.}
\newtheorem{theorem}{\bf Theorem}
\newtheorem{lemma}[theorem]{\bf Lemma}
\newtheorem{corollary}[theorem]{\bf Corollary}

\newtheorem{proposition}[theorem]{\bf Proposition}
\newtheorem{problem}{{\bf Problem}}

\newcommand{\epf}{\Box\vspace{0.15in}}

\newcommand{\gone}{{\sf G}_1}

\newcommand{\open}{\mathcal{O}}

\begin{document}
\maketitle
\begin{abstract} We obtain from the consistency of the existence of a measurable cardinal the consistency of ``small" upper bounds on the cardinality of a large class of Lindel\"of spaces whose singletons are ${\sf G}_{\delta}$ sets.
\end{abstract}
Call a topological space in which each singleton is a ${\sf G}_{\delta}$ set a \emph{points ${\sf G}_{\delta}$ space}. A.V. Arhangel'skii proved that any points ${\sf G}_{\delta}$ Lindel\"of space must have cardinality less than the least measurable cardinal and asked whether for ${\sf T}_2$ spaces this cardinality upper bound could be improved.  I. Juhasz constructed examples showing that for ${\sf T}_1$ spaces this upper bound is sharp. F.D. Tall, investigating Arhangel'skii's problem, defined the class of \emph{indestructibly} Lindel\"of spaces. A Lindel\"of space is \emph{indestructible} if it remains Lindel\"of after forcing with a countably closed forcing notion. He proved:
\begin{theorem}[F.D. Tall \cite{FDTGdelta}]\label{supercompact} If it is consistent that there is a supercompact cardinal, then it is consistent that $2^{\aleph_0} = \aleph_1$, and every points ${\sf G}_{\delta}$ indestructibly Lindel\"of space has cardinality $\le \aleph_1$.
\end{theorem}

In this paper we show that the hypothesis that there is a supercompact cardinal can be weakened to the hypothesis that there exists a measurable cardinal. Our technique permits flexibility on the cardinality of the continuum.

In Section 1 we review relevant information about ideals and the weakly precipitous ideal game. The relevance of the weakly precipitous ideal game to points ${\sf G}_{\delta}$ spaces is given in Lemma \ref{idealelements2}. In Section 2 we consider the indestructibly Lindel\"of spaces. A variation of the weakly precipitous ideal game is introduced. This variation is featured in the main result, Theorem \ref{TWOtacticandindL}: a cardinality restriction is imposed on the indestructibly Lindel\"of spaces with points ${\sf G}_{\delta}$. In Section 3 we give the consistency strength of the hypothesis used in Theorem \ref{TWOtacticandindL} and point out that mere existence of a precipitous ideal is insufficient to derive a cardinality bound on the indestructibly Lindel\"of points ${\sf G}_{\delta}$ spaces. In Section 4 we describe models of set theory in which the Continuum Hypothesis fails while there is a ``small" upper bound on the cardinality of points ${\sf G}_{\delta}$ indestructibly Lindel\"of spaces.

The notion of a Rothberger space appears in the paper. A space $X$ is a \emph{Rothberger space} if for each $\omega$-sequence of open covers of $X$ there is a sequence of open sets, then $n$-th belonging to the $n$-th cover, such that the terms of the latter sequence is an open cover of $X$. Rothberger spaces are indestructibly Lindel\"of (but not conversely). More details about Rothberger spaces in this context can be found in \cite{MSFT}.

\section*{Acknowledgement}

I would like to thank Professor F.D. Tall for very a stimulating correspondence regarding the problems treated in this paper, and for permission to include his argument (slightly adapted) in the proof of Theorem \ref{indlindsubsp}. I would also like to thank Professor Masaru Kada for important remarks about Theorem \ref{TWOtacticandindL}, and for permission to include his remarks (See remark (3) in the last section of the paper).

\section{Weakly precipitous ideals and points ${\sf G}_{\delta}$ spaces.}

For $\kappa$ an infinite cardinal, $\mathcal{P}(\kappa)$ denotes the powerset of $\kappa$. A set $J\subseteq\mathcal{P}(\kappa)$ is said to be a \emph{free ideal} on $\kappa$ if: (i) each finite subset of $\kappa$ is an element of $J$, (ii) $\kappa$ is not an element of $J$, (iii) the union of any two elements of $J$ is an element of $J$, and (iv) if $B\in J$ then $\mathcal{P}(B)\subset J$. For a free ideal $J$ the symbol $J^+$ denotes $\{A\in\mathcal{P}(\kappa):A\not\in J\}$.

Let $\lambda\le \kappa$ be a cardinal number. The free ideal $J$ on $\kappa$ is said to be $\lambda$-\emph{complete} if: For each $\mathcal{A}\subset J$, if $\vert\mathcal{A}\vert<\lambda$ then $\bigcup\mathcal{A}\in J$. A free ideal which is $\omega_1$-complete is said to be $\sigma$-complete. 

For a free ideal $J$ on $\kappa$ Galvin \emph{et al.} \cite{gjm} investigated the game $G(J)$ of length $\omega$, defined as follows: Two players, ONE and TWO, play an inning per finite ordinal $n$. In inning $n$ ONE first chooses $O_n\in J^+$. TWO responds with $T_n\in J^+$. The players obey the rule that for each $n$, $O_{n+1}\subset T_n\subset O_n$. TWO wins a play
\[
  O_1,\, T_1,\, O_2,\, T_2,\, \cdots, O_n,\, T_n,\, \cdots
\]
if $\bigcap_{n<\omega}T_n\neq \emptyset$; else, ONE wins. 

It is easy to see that if $J$ is not $\sigma$-complete, then ONE has a winning strategy in ${\sf G}(J)$. It was shown in Theorem 2 of \cite{gjm} that $J\subseteq \mathcal{P}(\kappa)$ is a \emph{weakly precipitous} ideal if, and only if, ONE has no winning strategy in the game $G(J)$. We shall take this characterization of weak precipitousness as the definition. An ideal $J$ on $\mathcal{P}(\kappa)$ is said to be \emph{precipitous} if it is weakly precipitous and $\kappa$-complete. This distinction was not made in the earlier literature such as \cite{gjm} and \cite{jmmp}. The $\kappa$-completeness requirement appears to have emerged in \cite{jp}, and the ``weakly precipitous" terminology for the $\sigma$-complete case seems to have been coined in \cite{kts}. 

This game is related as follows to spaces in which each point is ${\sf G}_{\delta}$:

\begin{lemma}\label{idealelements2} Let $\kappa$ be a cardinal such that there is a weakly precipitous ideal $J\subset\mathcal{P}(\kappa)$. Let $X\supseteq\kappa$ be a topological space in which each point is a ${\sf G}_{\delta}$. Then for each $x\in X$ and each $B\in J^+$ and each sequence $(U_n(x):n<\omega)$ of neighborhoods of $x$ with $\{x\}=\cap_{n<\omega}U_n(x)$, there is a $C \subseteq B$ with $C \in J^+$ and an $n$ such that $U_n(x) \cap C \in J$
\end{lemma}
{\bf Proof:} For each $x$ in $X$ fix a sequence $(U_n(x):n<\omega)$ of open neighborhoods such that for each $n$ we have $U_{n+1}(x)\subseteq U_n(x)$, and $\{x\} = \cap_{n<\omega}U_n(x)$. 

Suppose that contrary to the claim of the lemma, there is an $x\in X$ and a $B\in J^+$ such that for each $C \subset B$ with $C\in J^+$ and for each n, $U_n(x) \cap C \in J^+$. Fix $x$ and $B$. Define a strategy $\sigma$ for ONE of the game $G^1(J)$ as follows:
ONE's first move is $\sigma(X) = (U_1(x)\cap B)\setminus\{x\}$. When TWO responds with a $T_1\subseteq \sigma(X)$ and $T_1\in J^+$, ONE plays $\sigma(T_1) = (U_2(x)\cap T_1)\setminus\{x\}$. When TWO responds with $T_2\subseteq \sigma(T_1)$, ONE plays $\sigma(T_1,T_2) = (T_2 \cap U_3(x))\setminus\{x\}$, and so on. 

Observe that $\sigma$ is a legitimate strategy for ONE. But since ONE has no winning strategy in $G(J)$, consider a $\sigma$-play lost by ONE, say
\[
  O_1,\, T_1,\, O_2,\, T_2,\, \cdots
\]
Then $\cap_{n=1}^{\infty}T_n \neq \emptyset$, implying that $\cap_{n=1}^{\infty}(U_n(x)\setminus\{x\})\neq\emptyset$, a contradiction.
$\epf$

\section{The cardinality of points ${\sf G}_{\delta}$ indestructibly Lindel\"of spaces.}

For a space $X$ define the game $\gone^{\omega_1}(\open,\open)$ as follows: Players ONE and TWO play an inning for each $\gamma<\omega_1$. In inning $\gamma$ ONE first chooses an open cover $O_{\gamma}$ of $X$, and then TWO chooses $T_{\gamma}\in O_{\gamma}$. A play
\[
  O_0,\, T_0,\, \cdots,\, O_{\gamma},\, T_{\gamma},\, \cdots
\]
is won by TWO if $\{T_{\gamma}:\gamma<\omega_1\}$ is a cover of $X$. Else, ONE wins.

In \cite{MSFT} we proved the following characterization of being indestructibly Lindel\"of: 
\begin{theorem}[\cite{MSFT}, Theorem 1]\label{indlindchar} A Lindel\"of space $X$ is indestructibly Lindel\"of if, and only if, ONE has no winning strategy in the game $\gone^{\omega_1}(\open,\open)$.
\end{theorem}

Of several natural variations on $G(J)$ we now need the following one: The game $G^+(J)$ proceeds like $G(J)$, but TWO wins a play only when $\bigcap_{n<\omega}T_n\in J^+$; else, ONE wins.
Evidently, if TWO has a winning strategy in $G^+(J)$ then TWO has a winning strategy in $G(J)$. 
Similarly, if ONE has no winning strategy in $G^+(J)$, then ONE has no winning strategy in $G(J)$. A winning strategy in $G^+(J)$ for TWO which depends on only the most recent move of ONE is said to be a winning \emph{tactic}. 

\begin{theorem}\label{TWOtacticandindL} Assume there is a free, $\sigma$-complete ideal $J$ on $\kappa$ such that TWO has a winning tactic in $G^+(J)$. Then each points ${\sf G}_{\delta}$ indestructibly Lindel\"of space has cardinality less than $\kappa$. 
\end{theorem}
{\flushleft{\bf Proof:}} Let $X$ be a points ${\sf G}_{\delta}$ Lindel\"of space with $\vert X\vert\ge\kappa$. Let $Y$ be subset of $X$ of cardinality $\kappa$ and let $J\subset\mathcal{P}(Y)$ be a free ideal such that TWO has a winning tactic 
$\sigma$ in $G^+(J)$. We define a winning strategy $F$ for ONE of the game $\gone^{\omega_1}(\open,\open)$ and then cite Theorem 1 of \cite{MSFT}:

For each $x\in X$ fix a sequence of neighborhoods $(U_n(x):n<\infty)$ such that for $m<n$ we have $U_m(x)\supset U_n(x)$, and $\{x\}=\bigcap_{n<\omega}U_n(x)$. First, ONE does the following: For each $x\in X$: Choose $D_x\subseteq Y$ and $n$ so that $D_x\not\in J$, $U_{n}(x)\cap D_x \in J$ and set $C(x) = \sigma(D_x)$. Choose $n(C(x),x)<\omega$ such that $C(x)\cap U_{n(C(x),x)}\in J$. ONE's first move in $\gone^{\omega_1}(\open,\open)$ is 
\[
  F(\emptyset) = \{U_{n(C(x),x)}(x): x\in X\}. 
\]

When TWO chooses $T_0\in F(\emptyset)$, fix $x_0$ with $T_0 = U_{n(C(x_0),x_0)}(x_0)$. Define $C_0 = C(x_0)$, $D_0 = D_{x_0}$. Since $C_0\in J^+$ we choose by Lemma \ref{idealelements2} for each $x\in X$ a $D_{x_0,x}$ and an $n$ with:
\begin{enumerate}
  \item{$D_{x_0,x} \in J^+$ and $D_{x_0,x}\subset C(x_0)$ and}
  \item{$U_{n}(x)\cap D_{x_0,x}\in J$.} 
\end{enumerate}
  Put $C(x_0,x) = \sigma(D_{x_0,x})$, and choose $n(C(x_0,x),x)$ with $C(x_0,x)\cap U_{n(C(x_0,x),x)}(x) \in J$. 
 ONE plays 
\[
  F(T_0) = \{U_{n(C(x_0,x),x)}(x):x\in X\}.
\]
When TWO plays $T_1\in F(T_0)$, fix $x_1$ so that $T_1 = U_{n(C(x_0,x_1),x_1)}(x_1)$. Define $C_1 = C(x_0,x_1)$ and $D_1 = D_{x_0,x_1}$. Since $C_1 \in J^+$ we choose by Lemma \ref{idealelements2} for each $x\in X$ a $D_{x_0,x_1,x}$ and an $n$ with:
\begin{enumerate}
  \item{$D_{x_0,x_1,x}\in J^+$ and $D_{x_0,x_1,x}\subset C(x_0,x_1)$ and}
  \item{$U_{n}(x)\cap D_{x_0,x_1,x}\in J$.} 
\end{enumerate} 
  Put $C(x_0,x_1,x)=\sigma(D_{x_0,x_1,x})$ and choose $n(C(x_0,x_1,x),x)$  with $C(x_0,x_1,x)\cap U_{n(C(x_0,x_1,x),x)}(x)\in J$.
 ONE plays 
\[
  F(T_0,T_1) = \{U_{n(C(x_0,x_1,x),x)}(x):x\in X\},
\]
and so on.

At a limit stage $\alpha<\omega_1$ we have descending sequences $(C_{\gamma}:\gamma<\alpha)$ and $(D_{\gamma}:\gamma<\alpha)$ of elements of $J^+$ as well as a sequence $(x_{\gamma}:\gamma<\alpha)$ of elements of $X$ such that:
\begin{enumerate}
  \item{For each $\gamma$, $C_{\gamma} = C(x_{\nu}:\nu\le\gamma)$ and $D_{\gamma}=D_{(x_{\nu}:\nu\le\gamma)}$;}
  \item{For each $\gamma$, $D_{\gamma+1}\subset C_{\gamma} = \sigma(D_{\gamma})$}
  \item{$T_{\gamma} = U_{n(C_{\gamma},x_{\gamma})}(x_{\gamma})$.}
\end{enumerate} 
Since $\alpha$ is countable choose a cofinal subset $(\gamma_n:n<\omega)$ of ordinals increasing to $\alpha$. Then as for each $n$ we have $C_{\gamma_n} = \sigma(D_{\gamma_n})$ we see thst $(C_{\gamma_n}:n<\omega)$ are moves by TWO, using the winning tactic $\sigma$, in $G^+(J)$. Thus we have $\cap_{\gamma<\alpha}C_{\gamma} = \cap_{n<\omega}C_{\gamma_n}\in J^+$.

Then by Lemma \ref{idealelements2} choose for each $x\in X$ a $D_{(x_{\nu}:\nu\le \gamma)\frown (x)}$ and $n$ such that
\begin{enumerate}
  \item{ $D_{(x_{\nu}:\nu\le \gamma)\frown (x)}\in J^+$ and $D_{(x_{\nu}:\nu\le \gamma)\frown (x)}\subset \cap_{\gamma<\alpha}C_{\gamma}$ and}
  \item{ $D_{(x_{\nu}:\nu\le \gamma)\frown (x)}\cap U_n(x) \in J$.}
\end{enumerate} 
 Put 
\[
  C((x_{\nu}:\nu<\alpha)\frown (x))=\sigma(D_{((x_{\nu}:\nu<\alpha)\frown (x))})
\]
 and choose $n(C((x_{\nu}:\nu<\alpha)\frown x),x)$ such that: 
$C((x_{\nu}:\nu<\alpha)\frown x)\cap U_{n(C((x_{\nu}:\nu<\alpha)\frown x)}(x)\in J$.

Then ONE plays
\[
  F(T_{\gamma}:\gamma<\alpha) = \{U_{n(C((x_{\nu}:\nu<\alpha)\frown x),x)}(x):x\in X\}.
\]
This defines a strategy for ONE of the game $\gone^{\omega_1}(\open,\open)$. To see that it is winning, suppose that on the contrary there is an $F$-play won by TWO, say
\[
  O_0, T_0,\cdots, O_{\gamma}, T_{\gamma}, \cdots, \gamma<\omega_1,
\]
where $O_0 = F(\emptyset)$ and for each $\gamma>0$, $O_{\gamma} = F(T_{\beta}:\beta<\gamma)$. Since TWO wins $\gone^{\omega_1}(\open,\open)$, $X=\cup_{\gamma<\omega_1}T_{\gamma}$. But $X$ is Lindel\"of, and so we find a $\beta<\omega_1$ with $X = \cup_{\gamma<\beta}T_{\gamma}$. But then $C_{\alpha} = C(x_{\nu}:\nu<\alpha)$, $\alpha<\beta$ occurring in the definition of $F$ are in $J^+$ and satisfy for $\alpha<\beta$ that $C_{\beta}\subset C_{\alpha}$. It follows that for each $\gamma<\beta$ we have $T_{\gamma}\cap C_{\beta}\in J$, and as $J$ is $\sigma$-complete if follows that the $T_{\gamma}$ do not cover $C_{\beta}\subset X$, a contradiction. $\Box$

At one point in the above proof we made use of the fact that TWO has a winning tactic in $G^+(J)$. It may be the case that the conclusion of the theorem can be deduced from simply assuming that TWO has a winning strategy in $G^+(J)$\footnote{In fact, the conclusion of the theorem can be deduced from this formally weaker hypothesis. See Remark (3) at the end of the paper.}. One way to achieve this would be to show that if TWO has a winning strategy in $G^+(J)$, then TWO has a winning tactic in $G^+(J)$. I have not investigated this. 
\begin{problem} Let $J$ be a $\sigma$-complete free ideal on $\kappa$ such that TWO has a winning strategy in $G^+(J)$. Does it follow that TWO has a winning tactic in $G^+(J)$?
\end{problem}
Note that by Theorem 7 of \cite{GT}, if TWO has a winning strategy in $G^+(J)$, then TWO has a winning strategy $\sigma$ such that $T_1 = \sigma(O_1)$, and for each $n$, $T_{n+1} = \sigma(T_n,O_{n+1})$.
Here is another natural question which I have not explored:
\begin{problem} Let $J$ be a $\sigma$-complete free ideal on $\kappa$ such that TWO has a winning strategy in $G(J)$. Does it follow that TWO has a winning strategy in $G^+(J)$?
\end{problem}

\section{The hypothesis ``TWO has a winning tactic in $G^+(J)$".}

We now consider the strength of the hypothesis that TWO has a winning tactic in $G^+(J)$. First recall some concepts. An ideal $J$ is said to be $\lambda$-\emph{saturated} if whenever $\mathcal{B}\subset\mathcal{P}(\kappa)\setminus J$ is such that whenever $X\neq Y$ are elements of $\mathcal{B}$ then $X\cap Y\in J$, then $\vert\mathcal{B}\vert<\lambda$. We write $sat(J) = \min\{\lambda:J \mbox{ is $\lambda$-saturated}\}$. When $sat(J)$ is infinite it is regular and uncountable. Note that if $\lambda<\mu$ and if $J$ is $\lambda$-saturated then it is $\mu$-saturated.
It is well-known that every $\kappa^+$-saturated $\kappa$-complete ideal on $\kappa$ is precipitous (see Lemma 22.22 of \cite{jechme}). 

Next, let $J\subset\mathcal{P}(\kappa)$ be a $\sigma$-complete ideal and let $\lambda\le\kappa$ be an initial ordinal. For subsets $X$ and $Y$ of $\kappa$ write $X \equiv Y \mbox{ mod } J$ if the symmetric difference of $X$ and $Y$ is in $J$. Then $\mathcal{P}(\kappa)/J$ denotes the set of equivalence classes for this relation, and $\lbrack X\rbrack_J$ denotes the equivalence class of $X$.

A subset $D$ of the Boolean algebra $\mathcal{P}(\kappa)/ J$ is said to be \emph{dense} if there is for each $b\in\mathcal{P}(\kappa)/J$ a $d\in D$ with $d<b$. A dense set $D\subset\mathcal{P}(\kappa)/J$ is said to be $\lambda$-\emph{dense} if for each $\beta<\lambda$, for each $\beta$-sequence $b_0> b_1> \cdots > b_{\gamma} >\cdots$, $\gamma<\beta<\lambda$ of elements of $D$ there is a $d\in D$ such that for all $\gamma<\beta$, $d<b_{\gamma}$.

The \emph{Dense Ideal Hypothesis for $\kappa\ge\lambda$}, denoted {\bf DIH}($\kappa,\lambda$), is the statement: 
\begin{quote} There is a $\sigma$-complete free ideal $J\subset\mathcal{P}(\kappa)$ such that the Boolean algebra $\mathcal{P}(\kappa)/J$ has a $\lambda$-dense subset $D$.
\end{quote}
Note that if $\mu<\lambda$ then {\bf DIH}($\kappa,\lambda)\Rightarrow$ {\bf DIH}($\kappa,\mu)$. 

Consider the following five statements:
\begin{itemize}
  \item [I  ]{There exists a measurable cardinal}
  \item [II ]{There is an $\omega_1$ dense free ideal $J$ on an infinite set $S$.}
  \item [III]{There is a free ideal $J$ on a set $S$ such that TWO has a winning tactic in ${\sf G}^+(J)$.} 
  \item [IV  ]{There is a precipitous\footnote{Indeed, weakly precipitous works here.} ideal $J$ on an infinite set $S$.}
  \item [V ]{There is an $\kappa^+$-saturated $\kappa$-complete free ideal on a regular cardinal $\kappa$.}
\end{itemize} 

Then I$\Rightarrow$II (let $J$ be the dual ideal of the ultrafilter witnessing measurability), II $\Rightarrow$ III (see remarks (1), (3) and (4) on page 292 of \cite{gjm}), and evidently III $\Rightarrow$ IV. As already noted, V $\Rightarrow$ IV. 

In ZFC, for a statement P, let CON(P) denote ``P is consistent". Then we have CON(I) if, and only if, CON(IV), and CON(V) implies CON(I):
\begin{proposition}\label{measurableequi} The existence of a free ideal $J$ on an uncountable cardinal such that TWO has a winning tactic in $G^+(J)$ is equiconsistent with the existence of a measurable cardinal. 
\end{proposition}
{\bf Proof:} When TWO has a winning tactic in $G^+(J)$, then TWO has a winning strategy in $G(J)$, and thus ONE has no winning strategy in $G(J)$. It follows that $J$ is a weakly precipitous ideal. Jech \emph{et al.} \cite{jmmp} show that the existence of a weakly precipitous ideal is equiconsistent with the existence of a measurable cardinal. This shows that consistency of the existence of a free ideal $J$ on an uncountable set, such that TWO has a winning tactic in $G^+(J)$ implies the consistency of the existence of a measurable cardinal.

For the other direction: The argument in $\S$4 of \cite{gjm} can be adapted to show that if it is consistent that there is a measurable cardinal $\kappa$, then for any infinite cardinal $\lambda<\kappa$ it is consistent that $\mathbf{DIH}(\lambda^{++},\lambda^+)$ holds\footnote{The model in \cite{gjm} is obtained as follows: For $\kappa$ a measurable in the ground model, collapse all cardinals below $\kappa$ to $\aleph_1$ using the Levy collapse. One can show that with $\mu<\kappa$ an uncountable regular cardinal, collapsing all cardinals below $\kappa$ to $\mu$ produces a model of ${\mathbf{DIH}}(\mu^+,\mu)$, by verifying that Lemmas 1, 2 and 3 and the subsequent claims in \cite{gjm} apply \emph{mutatis mutandis}.}. A free ideal $J$ on $\lambda^{++}$ witnessing $\mathbf{DIH}(\lambda^{++},\lambda^+)$ is a free $\sigma$-complete ideal such that TWO has a winning tactic in $G^+(J)$. $\epf$

Property V is preserved when adding $\aleph_1$ or more Cohen reals. This follows from the following well-known fact stated as Lemma 22.32 in \cite{jechme}: 
\begin{lemma}\label{precipitouspreserve} Let $\mathbb{B}$ be a complete Boolean algebra, let $G$ be $V$-generic on $\mathbb{B}$ and let $\kappa$ be a regular uncountable cardinal. Assume that $sat(\mathbb{B})\le\lambda$ and $sat(\mathbb{B})<\kappa$. Then: If $J\subseteq\mathcal{P}(\kappa)$ is $\lambda$-saturated and $\kappa$-complete, then in $V\lbrack G\rbrack$ $J$ generates a $\lambda$-saturated $\kappa$-complete ideal on $\kappa$.
\end{lemma}

Now I. Juhasz has proven that for each infinite cardinal $\kappa$ less than the first measurable cardinal there is a points ${\sf G}_{\delta}$ Lindel\"of space $X$ with $\kappa<\vert X\vert$ (see \cite{FDTGdelta} for details). But adding $\aleph_1$ Cohen reals converts each such groundmodel space to a Rothberger space (and thus indestructibly Lindel\"of space) in the generic extension (see \cite{MSFT}). Thus if it is consistent that there is a $\mu^+$-saturated $\mu$-complete ideal on some regular cardinal $\mu$, then it is consistent that there is a (weakly) precipitous ideal on $\mu$, and yet there is an indestructibly Lindel\"of points ${\sf G}_{\delta}$ space of cardinality larger than $\mu$. 

\section{The continuum and the cardinality of points ${\sf G}_{\delta}$ indestructibly Lindel\"of spaces.}

The first consequence of the work above is that the hypothesis of the consistency of the existence of a supercompact cardinal in Theorem \ref{supercompact} can be reduced to the hypothesis of the consistency of the existence of a measurable cardinal:
\begin{corollary}\label{MS} If it is consistent that there is a measurable cardinal, then it is consistent that $2^{\aleph_0}=\aleph_1$ and all indestructible points ${\sf G}_{\delta}$ Lindel\"of spaces are of cardinality $\le \aleph_1$.
\end{corollary}

In what follows we demonstrate that the a bound on the cardinality of points ${\sf G}_{\delta}$ indestructibly Lindel\"of spaces does not have a strong influence on the cardinality of the real line.
\begin{corollary}\label{MS2} If it is consistent that there is a measurable cardinal $\kappa$, then for each regular cardinal $\aleph_{\alpha}$ with $\kappa>\aleph_{\alpha}>\aleph_0$ it is consistent that $2^{\aleph_0}=\aleph_{\alpha+1}$ and for all $\mu\le 2^{\aleph_0}$ there are indestructible points ${\sf G}_{\delta}$ Lindel\"of spaces and there are none of cardinality $> 2^{\aleph_0}$.
\end{corollary}
{\flushleft{\bf Proof:}} First raise the continuum to $\aleph_{\alpha+1}$ by adding reals. Next L\'evy collapse the measurable cardinal to $\aleph_{\alpha+2}$. In the resulting model $2^{\aleph_0}=\aleph_{\alpha+1}$ and ${\sf DIH}(\aleph_{\alpha+2},\aleph_{\alpha+1})$ holds. By Theorem 1 each indestructibly Lindel\"of space with points ${\sf G}_{\delta}$ has cardinality $\le \aleph_{\alpha+1}$ in this generic extension. Since each separable metric space is indestructibly Lindel\"of it follows that there is for each cardinal $\lambda\le \aleph_{\alpha+1}$ an indestructibly Lindel\"of points ${\sf G}_{\delta}$ space of cardinality $\lambda$. $\Box$

\begin{corollary}\label{MS3} If it is consistent that there is a measurable cardinal $\kappa$, then for each pair of regular cardinals $\aleph_{\alpha}<\aleph_{\beta}<\kappa$ with $\aleph_{\beta}^{\aleph_1} = \aleph_{\beta}$ it is consistent that $2^{\aleph_0}=\aleph_{\alpha}$ and $2^{\aleph_1}=\aleph_{\beta}$ and there is a points ${\sf G}_{\delta}$ indestructibly Lindel\"of space of cardinality $\aleph_{\beta}$, but there are no points ${\sf G}_{\delta}$ indestructible Lindel\"of spaces of cardinality $> 2^{\aleph_1}$.
\end{corollary}
{\flushleft{\bf Proof:}} We may assume the ground model is ${\mathbf L}\lbrack\mathcal{U}\rbrack$ where $\mathcal{U}$ is a normal ultrafilter witnessing measurability, and thus that GCH holds. First use Gorelic's cardinal- and cofinality- preserving forcing to raise $2^{\aleph_1}$ to $\aleph_{\beta}$ while maintaining CH. This gives a points ${\sf G}_{\delta}$ indestructibly Lindel\"of ${\sf T}_{3}$-space $X$ of cardinality $2^{\aleph_1}$. Then add $\aleph_{\alpha}$ Cohen reals to get $2^{\aleph_0}=\aleph_{\alpha}$. In this extension the space $X$ from the first step still is a points ${\sf G}_{\delta}$ indestructibly Lindel\"of ${\sf T}_{3}$-space since all these properties are preserved by Cohen reals \cite{MSFT}. The cardinal $\kappa$ is, in this generic extension, still measurable \cite{LS}. Finally, Levy collapse the measurable cardinal to $\aleph_{\beta+1}$. This forcing is countably closed (and more) and thus preserves indestructibly Lindel\"of spaces from the ground model. The resulting model is the one for the corollary.
$\Box$ 

\section{Regarding a problem of Hajnal and Juhasz.}

Hajnal and Juhasz asked if an uncountable ${\sf T}_2$-Lindel\"of space must contain a Lindel\"of subspace of cardinality $\aleph_1$. Baumgartner and Tall showed in \cite{BT} that there are ZFC examples of uncountable ${\sf T}_1$ Lindel\"of spaces with points ${\sf G}_{\delta}$ which have no Lindel\"of subspaces of cardinality $\aleph_1$. In \cite{KT}, Section 3, Koszmider and Tall showed that the answer to Hajnal and Juhasz's question is ``no". They also show that the existence of their example is independent of ZFC. Recall that a topological space is said to be a P-space if each ${\sf G}_{\delta}$-subset is open. It is known that Lindel\"of P-spaces are Rothberger spaces \cite{MSFT}.
They show in \cite{KT}, Theorem 4, that :
\begin{theorem}[Koszmider-Tall]\label{pspaces} The following is consistent relative to the consistency of ZFC: CH holds, $2^{\aleph_1}>\aleph_2$ and every ${\sf T}_2$ Lindel\"of P-space of cardinality $\aleph_2$ contains a convergent $\omega_1$-sequence (thus a Rothberger subspace of cardinality $\aleph_1$).
\end{theorem}

And then in Section 3 of \cite{KT} they obtain their (consistent) example:
\begin{theorem}[Koszmider-Tall]\label{Pspacenosubspace} It is consistent, relative to the consistency of ZFC, that CH holds and there is an uncountable ${\sf T}_3$-Lindel\"of P-space which has no Lindel\"of subspace of cardinality $\aleph_1$.
\end{theorem}

One may ask if the problem of Hajnal and Juhasz has a solution in certain subclasses of the class of Lindl\"of spaces. Koszmider and Tall's results show that even in the class of Rothberger spaces the Hajnal-Juhasz problem has answer ``no". In the class of Rothberger spaces with small character the following is known  \cite{MSFT}:
\begin{theorem}\label{supercompactrothb} If it is consistent that there is a supercompact cardinal, then it is consistent that $2^{\aleph_0} = \aleph_1$, and every uncountable Rothberger space of character $\le \aleph_1$ has a Rothberger subspace of cardinality $\aleph_1$. 
\end{theorem}

F.D. Tall communicated to me that the techniques of this paper can also be used to reduce the strength of the hypothesis in Theorem \ref{supercompactrothb} from supercompact to measurable. A small additional observation converts Tall's remark to the following. 
\begin{theorem}\label{indlindsubsp} Assume there is a free ideal $J$ on $\omega_2$ such that TWO has a winning tactic in $G^{+}(J)$. Then every indestructibly Lindel\"of space of cardinality larger than $\aleph_1$ and of character $\le\aleph_1$ has a Rothberger subspace of cardinality $\aleph_1$. 
\end{theorem}
{\bf Proof:} If an indestructibly Lindel\"of space has cardinality larger than $\aleph_1$ then Theorem \ref{TWOtacticandindL} implies it has  a point that is not ${\sf G}_{\delta}$. If a Lindel\"of space has character $\le\aleph_1$ and if some element is not a ${\sf G}_{\delta}$-point, then the space has a convergent $\omega_1$-sequence (Theorem 7 in \cite{BT}). Such a sequence together with its limit is a Rothberger subspace. $\Box$ 

\section{Remarks.}

{\flushleft{\bf (1)}} Lemma \ref{idealelements2} can be stated in greater generality that may be useful for other applications of these techniques:
\begin{lemma}\label{idealelements3} Let $\kappa$ be a cardinal such that there is a weakly precipitous ideal $J\subset\mathcal{P}(\kappa)$. Let $X\supseteq\kappa$ be a topological space and let $\mathcal{F}$ be a family of ${\sf G}_{\delta}$ subsets of $X$ such that $\mathcal{F}\subseteq J$. Then for each $F\in \mathcal{F}$ and each $B\in J^+$ and each sequence $(U_n(F):n<\omega)$ of open neighborhoods of $F$ with $F=\cap_{n<\omega}U_n(F)$, there is a $C \subseteq B$ with $C \in J^+$ and an $n$ such that $U_n(F) \cap C \in J$
\end{lemma}

{\flushleft{\bf (2)}} If in a ground model $V$ we have an ideal $J$ on an ordinal $\alpha$, then in generic extensions of $V$ let $J^*$ denote the ideal on $\alpha$ generated by $J$.  It is of interest to know which forcings increase $2^{\aleph_1}$ but preserve for example the statement: ``There is a $\sigma$-complete free ideal $J$ on $\omega_2$ such that TWO has a winning tactic in $G^+(J^*)$". 
This is not preserved by all $\omega_1$-complete $\omega_2$-chain condition partial orders: In \cite{IG} Gorelic shows that for each cardinal number $\kappa>\aleph_1$ it is consistent that CH holds, that $2^{\aleph_1}>\kappa$, and there is a ${\sf T}_3$ points ${\sf G}_{\delta}$-indestructibly Lindel\"of space $X$ of cardinality $2^{\aleph_1}$. Since the model in Section 4 of \cite{gjm} is a suitable ground model for Gorelic's construction, Theorem \ref{TWOtacticandindL} implies that in the model obtained by applying Gorelic's extension to the model from \cite{gjm}, there is no free ideal $J$ on $\omega_2$ such that TWO has a winning tactic in $G^+(J)$.

{\flushleft{\bf (3)}} After learning of the proof of Theorem \ref{TWOtacticandindL}, Masaru Kada informed me that in fact, by known results of Foreman and independently Veli\v{c}kovic, the hypothesis that TWO has a winning strategy in ${\sf G}^+(J)$ is sufficient to prove this theorem. We thank Kada for his kind permission to include the relevant remarks here. For notation and more information, see \cite{IY}: If TWO has a winning strategy in ${\sf G}^+(J)$, the the partially ordered set $(J^+,\subseteq)$ is $\omega+1$-strategically closed, and thus (see Corollary 3.2 in \cite{IY}) is strongly $\omega_1$-strategically closed. In the proof of Theorem \ref{TWOtacticandindL} use TWO's winning strategy in the game ${\sf G}^{I}_{\omega_1}(J^+)$ instead of a winning tactic in ${\sf G}^+(J)$. Note that also the existence of a free ideal $J$ for which TWO has a winning strategy in ${\sf G}^+(J)$ is equiconsistent with the existence of a measurable cardinal.

\end{document}